\documentclass[12pt]{article}

\usepackage{scicite}
\usepackage{todonotes}

\usepackage[CaptionAfterwards]{fltpage}

\usepackage[rightcaption]{sidecap}
\sidecaptionvpos{figure}{c}
\usepackage{epsfig}
\usepackage{graphicx}
\usepackage{amsmath}
\usepackage{amssymb}
\usepackage{wrapfig}
\usepackage{lipsum}
\usepackage{caption}
\usepackage{subcaption}
\usepackage{float}

\usepackage{times}
\usepackage{physics}
\usepackage{amsfonts}
\usepackage{xcolor}
\usepackage{makecell}
\usepackage{hyperref}

\definecolor{azure(colorwheel)}{rgb}{0.0, 0.5, 1.0}
\hypersetup{
    colorlinks=true,
    urlcolor=azure(colorwheel),
}

\newcommand{\Lagr}{\mathcal{L}}
\newcommand{\vecnorm}[1]{\left\lVert#1\right\rVert}

\usepackage{environ}

\newenvironment{sciabstract}{%
\begin{quote} \bf}
{\end{quote}}

\newcounter{lastnote}

\title{Discovering State Variables Hidden in Experimental Data}

\author
{Boyuan Chen$^{\ast}$, Kuang Huang, Sunand Raghupathi, \\
Ishaan Chandratreya, Qiang Du, Hod Lipson\\
\\
\href{https://www.cs.columbia.edu/~bchen/neural-state-variables}{neural-state-variables.com}\\
Columbia University\\
\\
\normalsize{$^\ast$To whom correspondence should be addressed; E-mail:  bchen@cs.columbia.edu.}
}

\date{}

\begin{document} 


\baselineskip18pt


\maketitle

\begin{sciabstract}
All physical laws are described as relationships between state variables that give a complete and non-redundant description of the relevant system dynamics. 
However, despite the prevalence of computing power and AI, the process of identifying the hidden state variables themselves has resisted automation. Most data-driven methods for modeling physical phenomena still assume that observed data streams already correspond to relevant state variables. A key challenge is to identify the possible sets of state variables from scratch, given only high-dimensional observational data. Here we propose a new principle for determining how many state variables an observed system is likely to have, and what these variables might be, directly from video streams. We demonstrate the effectiveness of this approach using video recordings of a variety of physical dynamical systems, ranging from elastic double pendulums to fire flames. Without any prior knowledge of the underlying physics, our algorithm discovers the intrinsic dimension of the observed dynamics and identifies candidate sets of state variables. We suggest that this approach could help catalyze the understanding, prediction and control of increasingly complex systems.
\end{sciabstract}

\clearpage

\subsubsection*{Main Text}

Mathematical relationships are known to describe nearly all physical laws in nature \cite{anderson1972more}, and these mathematical expressions are almost always formulated as relationships between physical state variables that describe the physical system. This suggests that before any natural law can be discovered, the relevant state variables must first be identified \cite{thompson2002nonlinear, hirsch2012differential}.

For example, it took civilizations millennia to formalize basic mechanical variables such as mass, momentum and acceleration. Only once these notions were formalized, could laws of mechanical motion be discovered. Laws of thermodynamics were discovered only after concepts such as temperature, pressure, energy and entropy were formalized. Laws of solid mechanics could only be discovered once variables such as stress and strain were formalized. Electromagnetism, fluid dynamics, quantum mechanics and so forth all required their own set of fundamental state variables to be defined, before they could be formalized into existence. Without the proper state variables, even a simple system may appear enigmatically complex.

\begin{figure*}[t!]
    \centering
    \includegraphics[width=0.7\linewidth]{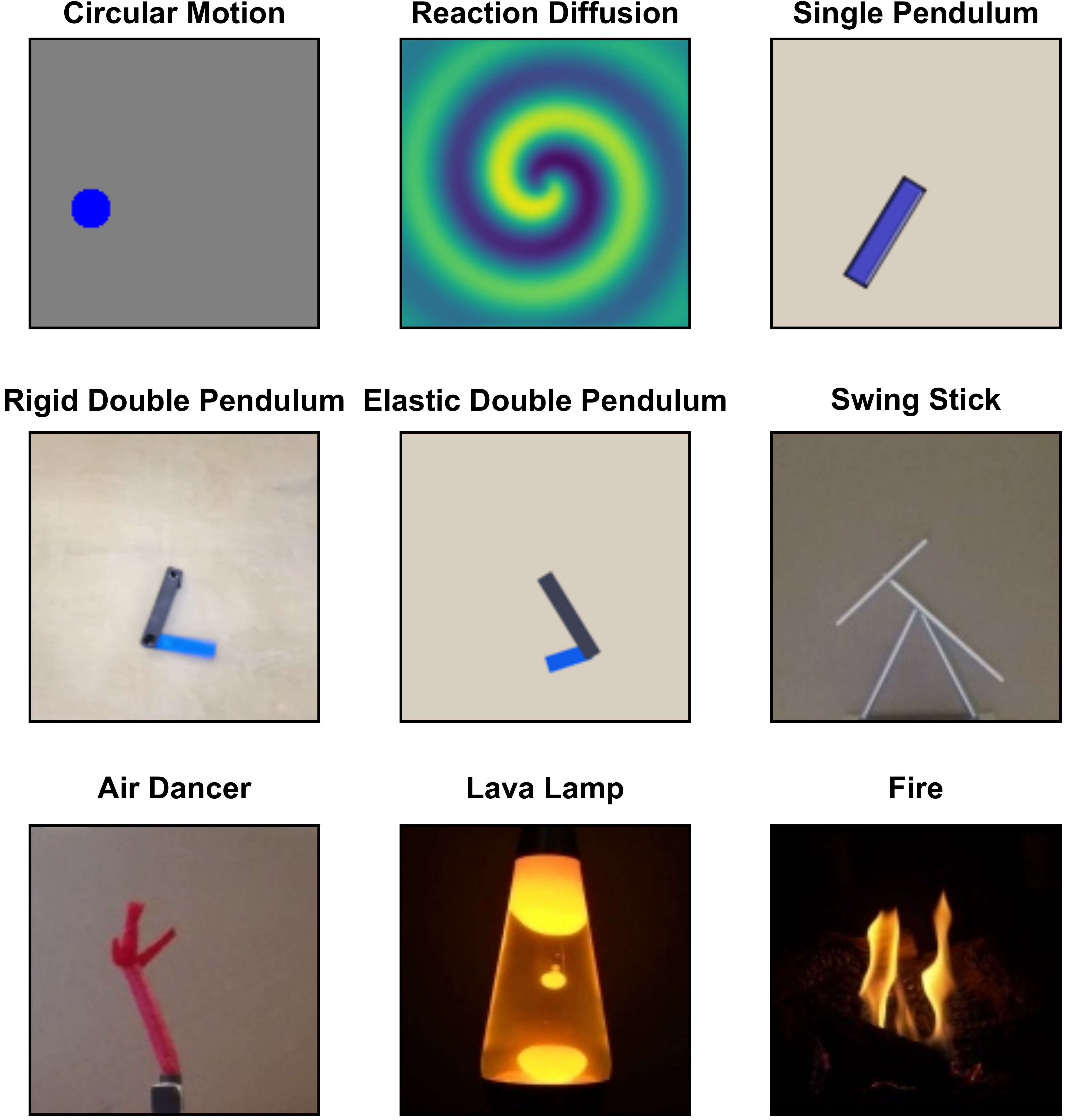}
    \caption{\textbf{What state variables shall we use to describe these dynamical systems?} Identifying state variables from raw observation data is a fundamental step to discover physical laws. The key challenge is to figure out what are the candidate variables, how the variables are dependent on each other, or how many variables will give a complete and non-redundant description of the system's states. Our work studies how to retrieve the possible set of state variables for data distributions non-linearly embedded in the ambient space.}
    \label{fig:teaser}
    \vspace{-12pt}
\end{figure*}

The set of state variables for modeling any system is not only hidden, but it is also not unique (Fig.~\ref{fig:teaser}). In fact, even for well-studied systems in classical mechanics, such as swinging pendulum, many sets of possible state variables exist. For the pendulum, the state variables are typically the angle of the arm $q_1 = \theta$ and the angular velocity of the arm $q_2 = \dot{\theta}$. The angle and angular velocity are convenient choices because they can be directly measured. However, alternative sets of state variables, such as kinetic and potential energies of the arm, could also be used. 

A key challenge, however, occurs when the system is new, unfamiliar or complex, and the relevant set of variables is unknown. 
Although there are various techniques such as Dynamic Mode Decomposition (DMD) and Singular Value Decomposition (SVD) \cite{kutz2016dynamic} developed to learn dynamical systems based on observations,  none of these methods has the ability to process a video of a pendulum, for example, and without any further knowledge, output the double pendulum's four state variables. Such an ability, if had, could help scientist gain insight into the physics underlying increasingly complex phenomena, when theory is not keeping pace with observations.

Data-analytics tools have impacted almost every aspect of scientific discovery \cite{evans2010machine, fortunato2018science}: Machines can measure, collect, store and analyze vast amounts of data. New machine learning techniques can create predictive models, find analytical equations and invariants, and even generate causal hypotheses along with new experiments to validate or refute these hypotheses \cite{king1996structure, waltz2009automating, king2009robot}. Yet, it remains challenging to develop automated systems to uncover the hidden state variables themselves. Finding such variables is still a laborious process requiring teams of human scientists toiling over decades.

The ability of human scientists to distill vast streams of raw observations into laws governing a concise set of relevant state variables has played a key role in many scientific discoveries. It is thus of great importance to have tools for automated scientific discovery that could help distill raw sensory perceptions into a compact set of state variables and their relationships.

Numerous machine learning tools have been demonstrated to model the dynamics of physical systems automatically, but most of them were already provided measurements of the relevant state variables in advance \cite{crutchfield1987equations, kevrekidis2003equation, yao2007modeling, bongard2007automated, rowley2009spectral, schmidt2009distilling, schmidt2011automated, sugihara2012detecting, ye2015equation, daniels2015automated, daniels2015efficient, benner2015survey, brunton2016discovering, rudy2017data, udrescu2020ai}\footnote{By 
State Variables, we refer to compact and complete sets of quantitative variables that fully describe the observed dynamical system evolving with time.}. For example, our own previous work on distilling natural laws \cite{schmidt2009distilling} assumed an input stream corresponding to state variables such as angle and angular velocity of a pendulum arm. Brunton et al. \cite{brunton2016discovering} required access to spatial coordinates and their derivatives for modeling a Lorenz system, Udrescu and Tegmark \cite{udrescu2020ai} combined neural networks with known physical properties to solve equations from the Feynman Lecture on Physics, given provided variables, Mrowca et al. \cite{mrowca2018flexible} required access to the position, velocity, mass, and material property of the particles that compose the objects being modeled and Champion et al. \cite{champion2019data} predefined possible basis functions to constrain the training of an autoencoder for observation reconstruction. 

The goal of this work is to find a way to overcome this key barrier to automated discovery – by explicitly identifying the intrinsic dimensionality of a system and the corresponding hidden state variables, purely from the visual information encoded in raw camera observations.
A key challenge in identifying state variables is that they are often hidden and might be hard to measure directly. An even greater challenging aspect of state variable identification is that there might be a large number of potential descriptive variables that are related to the varying state of the system, but are neither compact nor complete in their description of the system.

For example, a camera observing a swinging pendulum with an imaging resolution of $128 \times 128$ pixels in three color channels, will measure 49,152 variables per frame. Yet this enormous set of measurement, while intuitively descriptive, is neither compact nor complete: In fact, we know that the state of a swinging pendulum can be described fully by only two variables: its angle and angular velocity. Moreover these two state variables cannot be measured from a single video frame alone. In other words, a single frame, despite the large number of measurements, is insufficient to describe the full state of a pendulum.

The questions that we aim to answer are: Given a series of video frames of a swinging pendulum that contain the full and accurate motion trajectories, for example, is there a way to know that only two variables are required to describe its dynamics in full? And is there an automated process to reduce the vast deluge of irrelevant and superfluous pixel information into representations in terms of the two state variables? Naturally, we would like this process to work across a variety of physical systems and phenomena.

The starting point of our approach is to model the system dynamics directly from video representations via a neural network with bottleneck latent embeddings \cite{baldi1989neural, hinton1994autoencoders, masci2011stacked}. If the network is able to make accurate future predictions, the network should internally encapsulate a relationship connecting relevant current states with future states. Our paramount challenge is to distill and extract the hidden state variables from the network encoding.

Our key idea involves two major stages. First, after training the dynamics predictive neural network, we calculate the minimum number of independent variables needed to describe the dynamical systems, known as its {\em intrinsic dimension}, with geometric manifold learning algorithms. This initial stage produces accurate intrinsic dimension estimations on a variety of systems from the model's bottleneck latent embeddings which are already reduced by hundreds of times compared to the raw image space.

Armed with the intrinsic dimension obtained in the first stage, in the second stage, we design a latent reconstruction neural network to further identify the governing state variables with the exact dimension as the intrinsic dimension. We term these identified state variables {\em Neural State Variables}. Through both quantitative and qualitative experiments, we demonstrated that Neural State Variables can accurately capture the overall system dynamics.

Beyond our two-stage approach to reveal the system intrinsic dimension and the possible set of state variables, another major contribution of our work is to leverage the discovered Neural State Variables both as an intermediate representation and evaluation metric for stable long-term future predictions of system behaviors. Due to the special reduced-dimension property of Neural State Variables, they can provide very stable long-term predictions, while higher dimensional auto-encoders often yield blurred or plain background predictions if iterated just a few steps into the future. 

Finally, we present a hybrid prediction scheme which achieves both accurate and stable long-term predictions. Furthermore, we derive a quantitative evaluation metric for long-term prediction stability with Neural State Variables by approximating the true system dynamics using the most compact latent space. 

We also demonstrate that Neural State Variables can offer a robust representation space for modeling system dynamics under various visual perturbations.

\subsubsection*{Modeling Dynamical Systems from Videos}

The dynamics of a physical system defines the rule that governs how the current system states will evolve into their successive states in the future. Mathematically, provided the ambient space $\mathcal{X}$ and the state space $\mathcal{S}\subset\mathcal{X}$, one can formulate the dynamical system as
\begin{align}\label{eq:dynamic_system}
    \vb*{X}_{t+dt}=F(\vb*{X}_t), \quad t=0,dt,2dt,3dt,\dots,
\end{align}
where $\vb*{X}_t\in\mathcal{S}$ is the system's current state at time $t$, $dt$ is the discrete time increment. $F:\mathcal{S}\to\mathcal{S}$ describes the state evolution from $\vb*{X}_t$ to the system's successive state $\vb*{X}_{t+dt}$ at time $t+dt$. Throughout this paper we will consider the system as discrete in time. Any continuous in time dynamical system can also be discretized to formulation \eqref{eq:dynamic_system} with an appropriate sampling interval $dt$.

The first step towards modeling a dynamical system is to choose the representation of system states. Previous studies assume that the states are given as the direct measurements of a set of pre-defined state variables, such as the position and velocity of a rigid body object. However, defining which state variables to measure requires expert prior knowledge of the system. For an unfamiliar physical system, we do not know in advance what quantities to measure. Moreover, most state variables are not directly measurable, as they correspond to properties that are physically unobservable in a non-intrusive manner or cannot be uniquely determined without prior knowledge. One example is the measurement of the position of the pendulum arms is very challenging without installing tracking markers which can change the system dynamics.

In this work, we chose video frames as the state representation. Using the notations above, $\mathcal{X}$ is the high dimensional image space. This choice comes with several advantages. First, video recordings do not require prior knowledge of the inner working processes of the observed dynamical system. Second, video cameras collect a rich stream of physics signals, without requiring expensive and specialized equipment. If we can apply our method to data collected by video cameras, then this approach could potentially operate with other types of sensor arrays.

Our goal is to learn a the most compact space that implicitly captures the entire system dynamics, using only high dimensional visual data as input. To achieve this, we formulate a self-supervised learning problem to leverage the natural supervision from future video streams. Our model is based on an auto-encoder neural network to map the high-dimensional visual observations to a relative low-dimensional embedding which will then be projected to a high-dimensional image space again to predict the future video frames. If the model can successfully produce accurate future predictions, the low-dimensional bottleneck has to capture sufficient information on the system dynamics.

\begin{figure*}[t!]
    \centering
    \includegraphics[width=0.9\linewidth]{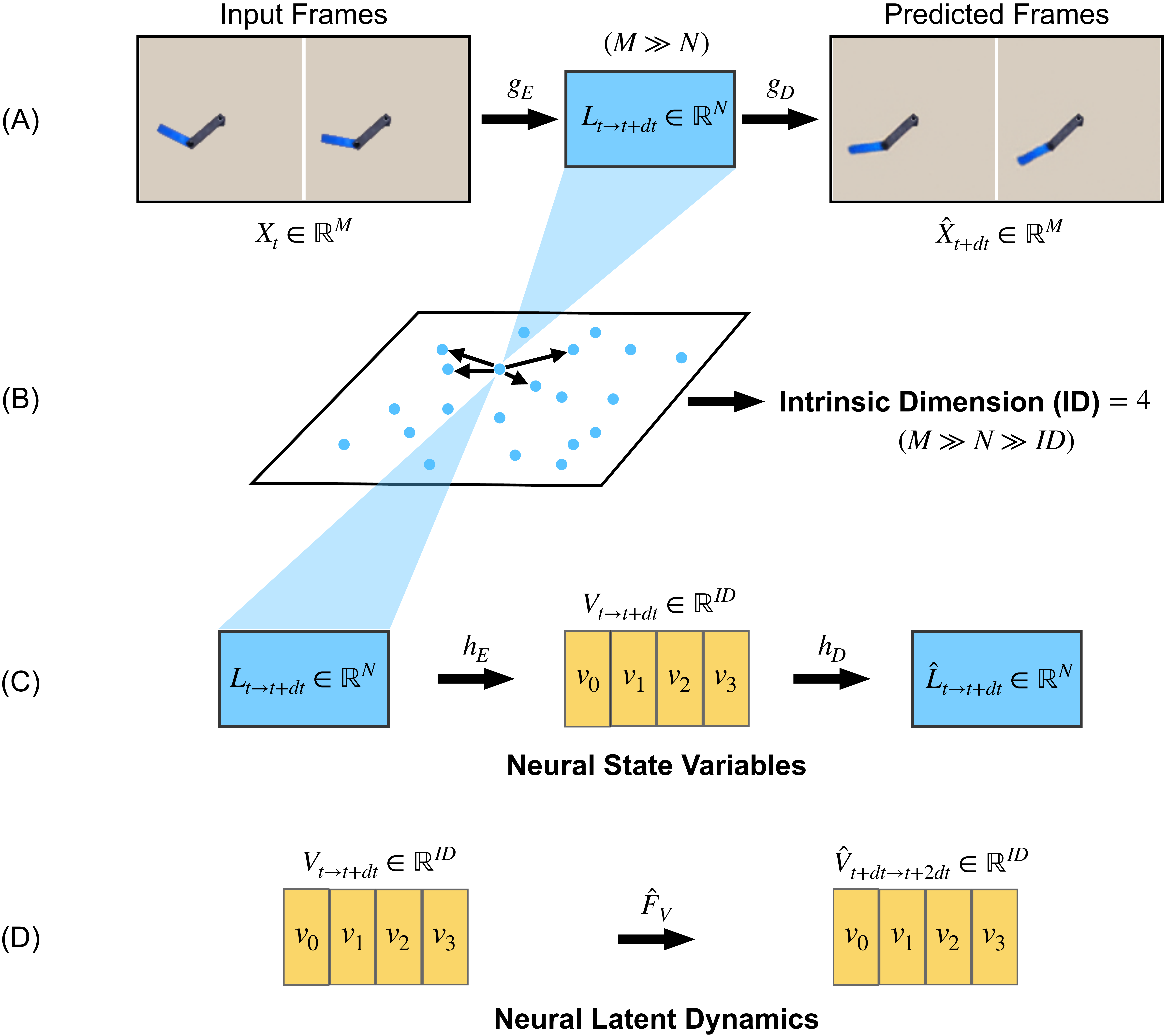}
    \caption{\textbf{Two-stage modeling of dynamical systems.} (A) and (B) \textbf{First stage: intrinsic dimension estimation.} We first modelled the dynamical systems via the evolution from $\vb*{X_t}$ to $\vb*{X_{t+dt}}$ with a fully convolutional encoder-decoder network directly from video observations. The dimension of the latent vectors $\vb*{L_{t \to t+dt}}$ is often much lower than the dimension of the input vectors, but much higher than the intrinsic dimension of the system. To identify the intrinsic dimension of the system, we applied geometric manifold learning algorithms on this set of relative high dimension latent vectors. (C) \textbf{Second stage: discover Neural State Variables.} We applied another encoder-decoder network on top of the above latent vectors to automatically determine the Neural State Variables by limiting the latent dimension of this network with the identified intrinsic dimension. Our two-step approach can produce Neural State Variables with the exact dimension of the system intrinsic dimension. (D) Once we determine the Neural State Variables, we can leverage the system dynamics in the space of Neural State Variables as an indicator of dynamics stability. Therefore, we learn a neural latent dynamics to predict the the Neural State Variables at the next time step from the current Neural State Variables.}
    \label{fig:overview}
\end{figure*}


Formally, our framework comprises five major components as shown in Fig.~\ref{fig:overview}(A): A pair of input image frames $\vb*{X_t}$, an encoder network $g_{\textrm{E}}$, a latent embedding vector $\vb*{L_{t \to t+dt}}$, a decoder network $g_{\textrm{D}}$, and a pair of output future frames $\vb*{X_{t+dt}}$. For the dynamical systems studied in our paper, both the input and output image pairs are two consecutive frames with the dimension $128 \times 128 \times 3$ RGB channels. The pairs of frames are concatenated to form single input and output image. The encoder $g_{\textrm{E}}$ and decoder $g_{\textrm{D}}$ are fully convolutional networks where the decoder also comes with residual connections and multi-scale predictions to enhance higher resolution predictions. The network first outputs $\vb*{L}_{t\to t+dt} = g_{\textrm{E}}(\vb*{X_t})$ and then generates the future frames $\vb*{\hat{X}_{t+dt}}$. 


To train the encode and decoder networks, we use simple L2 loss function without other constraints: 
$$\Lagr = \mathbb{E}_{\vb*{X}}\left[ {\lVert g_{\textrm{D}}(g_{\textrm{E}}(\vb*{X_t})) - \vb*{X_{t+dt}} \rVert_2^2}\right].$$ 
The learned mapping $\hat{F}=g_{\textrm{D}}\circ g_{\textrm{E}}$ provides a numerical approximation of the system's evolution mapping $F$ through the latent embedding.

One critical but largely ignored design decision is the dimension LD of the latent embedding $\vb*{L} \in \mathbb{R}^{\textrm{LD}}$. In Machine Learning, LD is often treated as a hyperparameter selected using an ``educated guess'' because it is not immediately clear what the best value of LD should be. However, this dimensionality is especially important for physical dynamics modeling. When LD is large, the latent embedding can hold large numbers of useful bits of information about the system dynamics. However large embedding vectors overfit the data are have limited capacity for longer range prediction. More importantly, large latent spaces hide and obfuscate the compact set of state variables we are after. When LD is too small, the network may under-fit the data. 

Therefore, we aim to come as close as possible to the exact number of state variables. In the next section, we will carefully analyze these issues and introduce our novel solutions to systematically discover the most compact space of the embedding vector without imposing optimization challenges.

To study the generality of the proposed approach, we compiled a dataset comprising videos of nine physical dynamical systems from various experimental domains (Fig.~\ref{fig:teaser}), ranging from simple periodic motion (circular motion, single pendulum), chaotic kinematics (rigid double pendulum, elastic double pendulum, swing stick), nonlinear wave (reaction-diffusion system), multi-phase flow (lava lamp), to aeroelasticity (air dancer) and combustion (flame dynamics). Each dataset comes with raw video recordings of the dynamical systems. For certain dataset where measurements of known physical quantities are easy to obtain, we also provide ground-truth estimations for evaluation purposes only. We include full details of the dataset in the Appendix.

\begin{FPfigure}
    \centering
    \includegraphics[width=\linewidth]{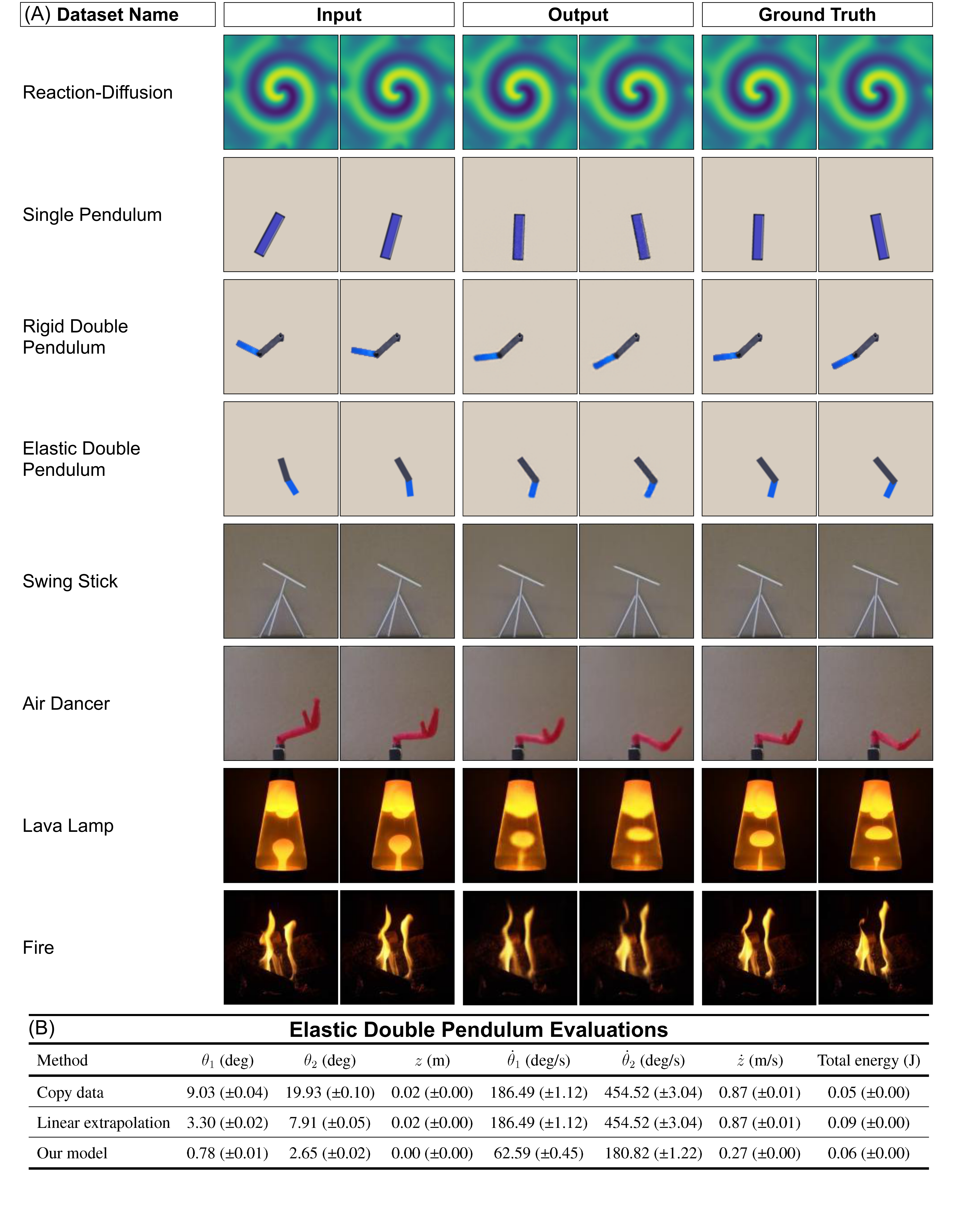}
    \caption{\textbf{Prediction visualizations and physics evaluations} (A) Visualizations of our basic prediction results. (B) For systems where physical variables happen to be available, we obtained the physical variables from both the predicted frames and the ground truth frames. We then performed physics evaluations on these systems. We show the results of elastic double pendulum here and include results for other systems in the Appendix. The elastic double pendulum dataset has 60 fps. Our prediction model outperforms both the copy data and linear extrapolation baselines suggesting that our model captures nontrivial understanding of the system's second order dynamics.}
    \label{fig:one-step}
    \vspace{-12pt}
\end{FPfigure}

In Fig.~\ref{fig:one-step}(A), we show comparisons between the predicted video frames and the ground-truth recordings. Our model was able to produce accurate video predictions. Our model also substantially outperforms linear extrapolation and copying input data baselines in Fig.~\ref{fig:one-step}(B). For dataset with ground-truth physical quantities such as elastic double pendulum, our system was able to predict the physical variables accurately compared to the ground truth. Overall, the evaluation results suggest that the model successfully captured a nontrivial understanding of the system dynamics.

\subsubsection*{Intrinsic Dimension Estimation}

Intrinsic Dimension (ID) has served as a fundamental concept in many advances in physical sciences. In general, the intrinsic dimension refers to the minimum number of independent variables needed to fully describe the state of a dynamical system. The intrinsic dimension is independent of specific representations of the system or choice of a particular set of state variables. In a more quantitative way, the intrinsic dimension could be equivalently defined as the topological dimension of the state space $\mathcal{S}$ as a manifold in the ambient space $\mathcal{X}$ \cite{bishop1995neural,camastra2016intrinsic,campadelli2015intrinsic}.

A common assumption when analyzing a physical system is that the intrinsic dimension is known a priory. An even stronger assumption is that the corresponding state variables themselves are given. Yet these assumptions do not hold for unknown or partially known systems. In order to uncover the underlying dynamics of a wide range of systems and make future predictions of their future behaviors, we need to automatically identify the intrinsic dimension of the systems and extract the corresponding state variables from observed data, which is often high dimensional and noisy.

\begin{figure*}[t!]
    \centering
    \includegraphics[width=0.9\linewidth]{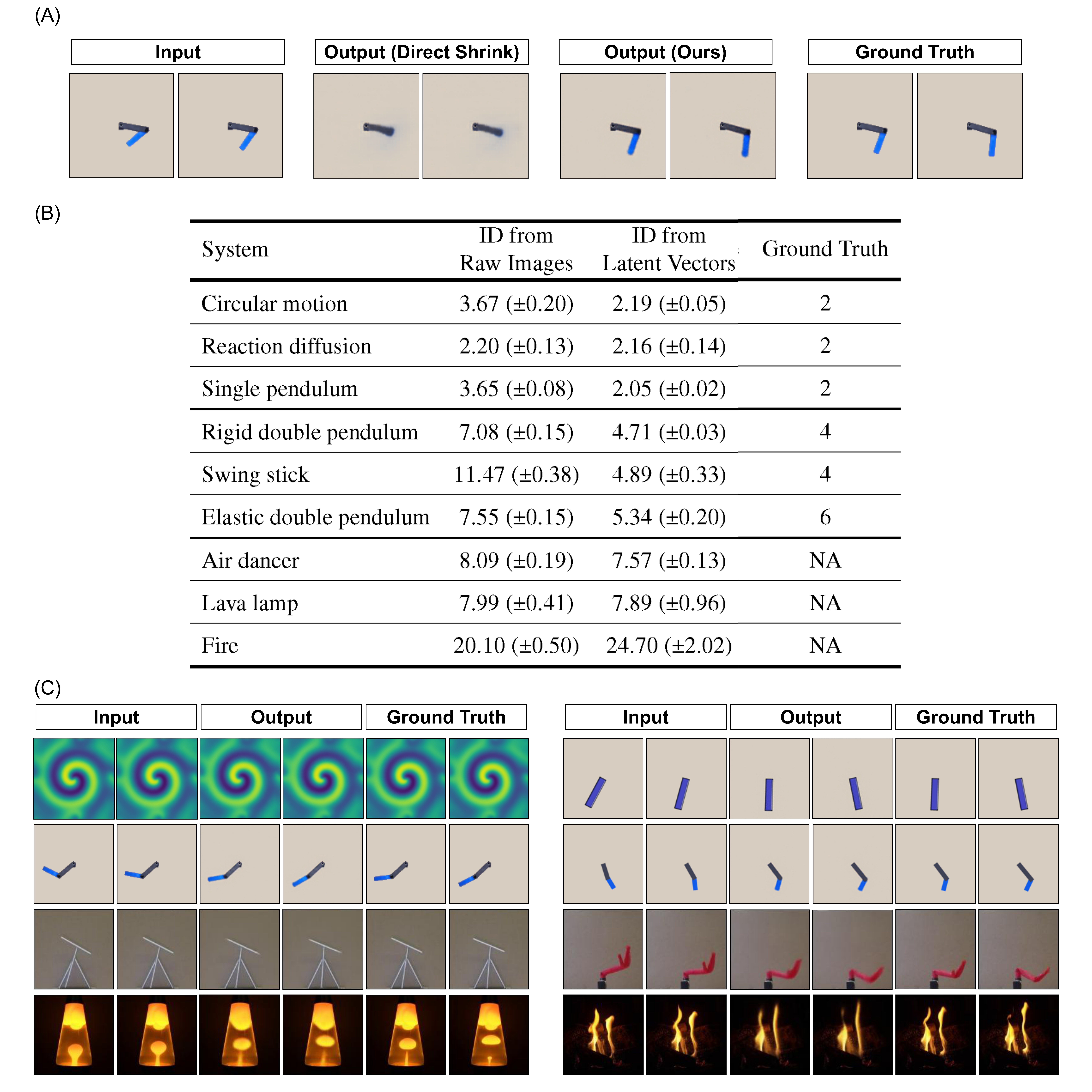}
    \caption{\textbf{Intrinsic Dimension (ID) and Neural State Variables} (A) Keep reducing the size of the latent embedding on the original auto-encoder to find ID is not feasible due to optimization difficulties. The network could not converge to a satisfactory solution. With our two-stage method to retrieve the system intrinsic dimension and further discovered Neural State Variables, we can bypass this 
    limitation to produce accurate future predictions. (B) Our method estimates ID without prior knowledge about the systems' state variables. The estimated ID value are rounded to the nearest even integer as position and velocity variables are in pairs. For systems with known IDs, our calculations give accurate results. For unknown systems, the ranking of the ID also makes sense. Our method outperforms direct estimations from raw images. (C) More results on one-step prediction with our discovered Neural State Variables.}
    \label{fig:intrinsic-dimension}
    \vspace{-12pt}
\end{figure*}

A naive approach using an auto-encoder predictive framework is to keep reducing the size of the latent embedding vector through trial and error until the output is no longer valid. However, this approach does not yield satisfactory results because the output deteriorates long before the minimal set of state variables is reached. As shown in Fig.~\ref{fig:intrinsic-dimension}(A), the model predictions broke down when we directly shrank the size of the latent space to the intrinsic dimension.

Inspired by traditional manifold learning methods that utilize geometric structures of the embedding vectors (such as their nearest distances), we propose a solution that can automatically discover the intrinsic dimension of a dynamical system from the latent vectors. Our approach only needs a one-time network training step. Specifically, we applied the Levina-Bickel’s algorithm \cite{levina2005maximum} on the latent embedding space. The algorithm considers latent vectors $\{\vb*{L}^{(1)},\vb*{L}^{(2)},\cdots,\vb*{L}^{(N)}\}$ collected from the trained dynamics predictive model as $N$ data points on a manifold of dimension ID in the latent embedding space. A key geometric observation is that the number of data points within distance $r$ from any given data point $\vb*{L}^{(i)}$ is proportional to $r^{\textrm{ID}}$ when $r$ is small. Based on the observation, the Levina-Bickel’s algorithm derives the local ID estimator near $\vb*{L}^{(i)}$ as $\frac{1}{k-2}\sum_{j=1}^{k-1}\log\frac{T_k(\vb*{L}^{(i)})}{T_j(\vb*{L}^{(i)})}$, where $T_k(\vb*{L}^{(i)})$ is the Euclidean distance between $\vb*{L}^{(i)}$ and its $k^{\textrm{th}}$ nearest neighbor in $\{\vb*{L}^{(1)},\vb*{L}^{(2)},\cdots,\vb*{L}^{(N)}\}$.
The global ID estimator is then calculated as:
\[ \textrm{ID}_{\textrm{L-B}} = \frac1N\sum_{i=1}^N\frac{1}{k-2}\sum_{j=1}^{k-1}\log\frac{T_k(\vb*{L}^{(i)})}{T_j(\vb*{L}^{(i)})}.\]

Fig.~\ref{fig:intrinsic-dimension}(B) shows the estimations across all the systems in our holdout dataset along with baseline comparisons from raw image observations and partial ground-truths. Our method demonstrates highly accurate estimations of the intrinsic dimension of all known systems. Although we cannot account for the ground-truth intrinsic dimension of other systems, we do see that our experiments presented a reasonable and intuitive relative ranking among all listed systems.

We also compared the performance of Levina-Bickel’s algorithm with other popular intrinsic dimensionality estimation algorithms including MiND\_ML, MiND\_KL, Hein, and CD \cite{campadelli2015intrinsic,rozza2012novel,ceruti2014danco,hein2005intrinsic,grassberger2004measuring} by following the original implementations \cite{hein2005intrinsic,ceruti2014danco}. We present full evaluations in the Appendix. Though all the algorithms demonstrated promising performance, we found that the Levina-Bickel algorithm gives the most robust and reliable estimation.

\subsubsection*{Neural State Variables}

As we have discussed above, the minimum set of independent state variables $\vb*{V}$ used to describe the dynamical system has the dimension known as the intrinsic dimension, namely $\vb*{V} \in \mathbb{R}^{\textrm{ID}}$. To simplify the terminology, we refer them as \textbf{State Variables} directly throughout the rest of this paper. 


Now that we have identified the number of state variables, we need to find the actual variables themselves (bearing in mind that the set is not unique). We propose a two-stage framework to retrieve possible State Variables from raw video data with both dynamics predictive and latent reconstruction neural networks. We term our subset of State Variables as \textbf{Neural State Variables}. Hence, with Neural State Variables $\vb*{V}$, the dynamical system can be expressed as the evolution of the trajectory $\{\vb*{V}_{0\to dt}$, $\vb*{V}_{dt \to 2dt}$, $... \}$.

Our key idea is to break down the identification process into two stages as illustrated in Fig.~\ref{fig:overview}(A)-(C). The first stage is to train a dynamics predictive model and identify the intrinsic dimension $\textrm{ID}$ as indicated in the previous section. This stage yields a relative low-dimensional latent embedding $\vb*{L} \in \mathbb{R}^{\textrm{LD}}$ where $\textrm{LD}$ is still much larger than $\textrm{ID}$. As shown in Fig.~\ref{fig:one-step}, the network can converge to output accurate future frames, indicating that the latent embedding captures sufficient information about the complete system dynamics.

The second stage operates directly on the latent embedding to further distill the Neural State Variables. Directly reducing the latent vector using the first step did not converge (Fig.~\ref{fig:intrinsic-dimension}(A)). Therefore, we trained a second auto-encoder network that takes in the pre-trained latent embedding and outputs the reconstruction of the input. The special property of this network is that the size of the latent embedding equals to the intrinsic dimension obtained from the first step. With a minimum reconstruction error, we can identify this latent embedding vectors as the Neural State Variables. Specifically, the network can be expressed as follows: $\vb*{V}_{t\to t+dt} = h_{\textrm{E}}(\vb*{L}_{t\to t+dt})$ and $\hat{\vb*{L}}_{t\to t+dt} = h_{\textrm{D}} (\vb*{V}_{t\to t+dt})$ where $h_{\textrm{E}}$ and $h_{\textrm{D}}$ refers to the encoder and decoder network of the latent reconstruction model. We train the latent reconstruction model with the L2 loss: $\Lagr = \mathbb{E}_{\vb*{L}}\left[ {\lVert h_{\textrm{D}}(h_{\textrm{E}}(\vb*{L_{t\to t+dt}})) - \vb*{L_{t\to t+dt}} \rVert_2^2}\right]$.

Overall, our two-stage method bypasses the optimization challenges and avoides the risk of under-fitting the observed data. In Fig.~\ref{fig:intrinsic-dimension}(C), we qualitatively demonstrate the effectiveness of our approach. For all the systems in our dataset, our framework is able to predict accurate future frames from super compact variables with dimension $\textrm{ID}$ (e.g., $\textrm{ID}=4$ for rigid double pendulum and $\textrm{ID}=6$ for elastic double pendulum). In the next section, we will provide quantitative evaluations and demonstrate its key usage to predict the long-term evolution of these dynamical systems.  

\subsubsection*{Neural State Variables for Stable Long-Term Prediction}

Forecasting the long-term future behaviors of unknown physical systems by learning to model their dynamics is critical for numerous real-world tasks. With a dynamics predictive model giving the one-step prediction, we can perform model rollout to feed each step's prediction as the input to predict the next state. However, there are two main challenges to obtaining satisfactory long-term predictions:

\begin{itemize}
    \item \textbf{One-Step Prediction Accuracy} The learned dynamics may not be accurate since prediction errors are iteratively introduced at every prediction step. This issue mainly attributes to the one-step prediction accuracy.
    \item \textbf{Long-term Prediction Stability} Due to error accumulation, the predicted sequences may not be able to maintain the ground truth state space: one repeated observation from past studies is that the long-term predicted sequences become blur, heavily distorted, or plain background within only a few rollouts. We also observed similar phenomena in our experiments as shown in Fig.~\ref{fig:intrinsic-dimension}(A). This is a very important issue to resolve because if objects deform or entirely disappear without following the system dynamics, it would be impossible to follow the system evolution faithfully. Here, we first define this phenomena as long-term prediction stability. We then present quantitative analysis for various prediction schemes. Finally, we will present our solution to approach the stability challenge with Neural State Variables.
\end{itemize}

Long-term prediction stability refers to the deviation between the predicted sequences generated from the learned dynamics and the ground truth state space governed by the system dynamics. Given a metric $M_{\mathcal{S}}(\cdot)$ that measures the deviation from a predicted state to the true state space $\mathcal{S}$, and a prediction sequence $\{\hat{\vb*{X}}_0,\hat{\vb*{X}}_{dt},\cdots\}$ from any initial state $\hat{\vb*{X}}_0$, we can quantify the stability of a prediction scheme as the growth rate of $M_{\mathcal{S}}(\hat{\vb*{X}}_t)$ as a function of $t$.

One challenge is to define at what point is the predicted image so degraded that it does not count as a prediction at all. We define an image quality test as follows (used only for evaluation): For systems for which we have prior knowledge about their conventional state variables, and we can extract these physical variables from the corresponding videos through classic computer vision techniques (e.g., color and contour extraction), $M_{\mathcal{S}}^{\textrm{phys}}(\cdot)$ can be readily defined as a binary value indicating whether the same set of physical variables can still be distilled from a predicted state $\hat{\vb*{X}}$ as its corresponding ground truth state. Consequently, if the predicted frame is heavily blurred or distorted, we will not be able to distill meaningful physical variables. Thus, $M_{\mathcal{S}}^{\textrm{phys}}(\hat{\vb*{X}})$ will be one. Otherwise $M_{\mathcal{S}}^{\textrm{phys}}(\hat{\vb*{X}})$ will be zero.

Moreover, to more generally capture the long-term predictive stability of various prediction schemes, $M_{\mathcal{S}}^{\textrm{phys}}(\cdot)$ should be evaluated on prediction sequences with multiple initial states. Therefore, we further define $M_{\mathcal{S}}^{\textrm{phys}}(\cdot)$ as the Reject Ratio to indicate how many predicted frames at each time step from different initial states will fail to pass the physical variables extraction test.

With the above test, we can quantitatively compare the stability of various long term prediction schemes. These schemes are based on iterative model rollouts but they differ in the size of intermediate variables.
When the model rollouts are through high dimensional latent vectors (8192 variables or 64 variables), the iterative scheme is given by $ \hat{\vb*{X}}_{t+dt}=g_{\textrm{D}}\circ g_{\textrm{E}}(\hat{\vb*{X}}_t),\quad t=0,dt,\dots$, where $g_{\textrm{E}}$ and $g_{\textrm{D}}$ represent the first auto-encoder that transforms input frames to the predicted frames via latent embeddings.
When the model rollouts are through Neural State Variables, the iterative scheme is given by $ \hat{\vb*{X}}_{t+dt}=g_{\textrm{D}}\circ h_{\textrm{D}}\circ h_{\textrm{E}} \circ  g_{\textrm{E}}(\hat{\vb*{X}}_t),\quad t=0,dt,\dots$, where $h_{\textrm{E}}$ and $h_{\textrm{D}}$ represent the second auto-encoder that reconstructs latent embeddings via Neural State Variables. The original latent embeddings are computed from input frames with $g_{\textrm{E}}$, and the reconstructed latent embeddings will be sent to $g_{\textrm{D}}$ to produce the final predicted frames.

\begin{figure*}[t!]
    \centering
    \includegraphics[width=0.9\linewidth]{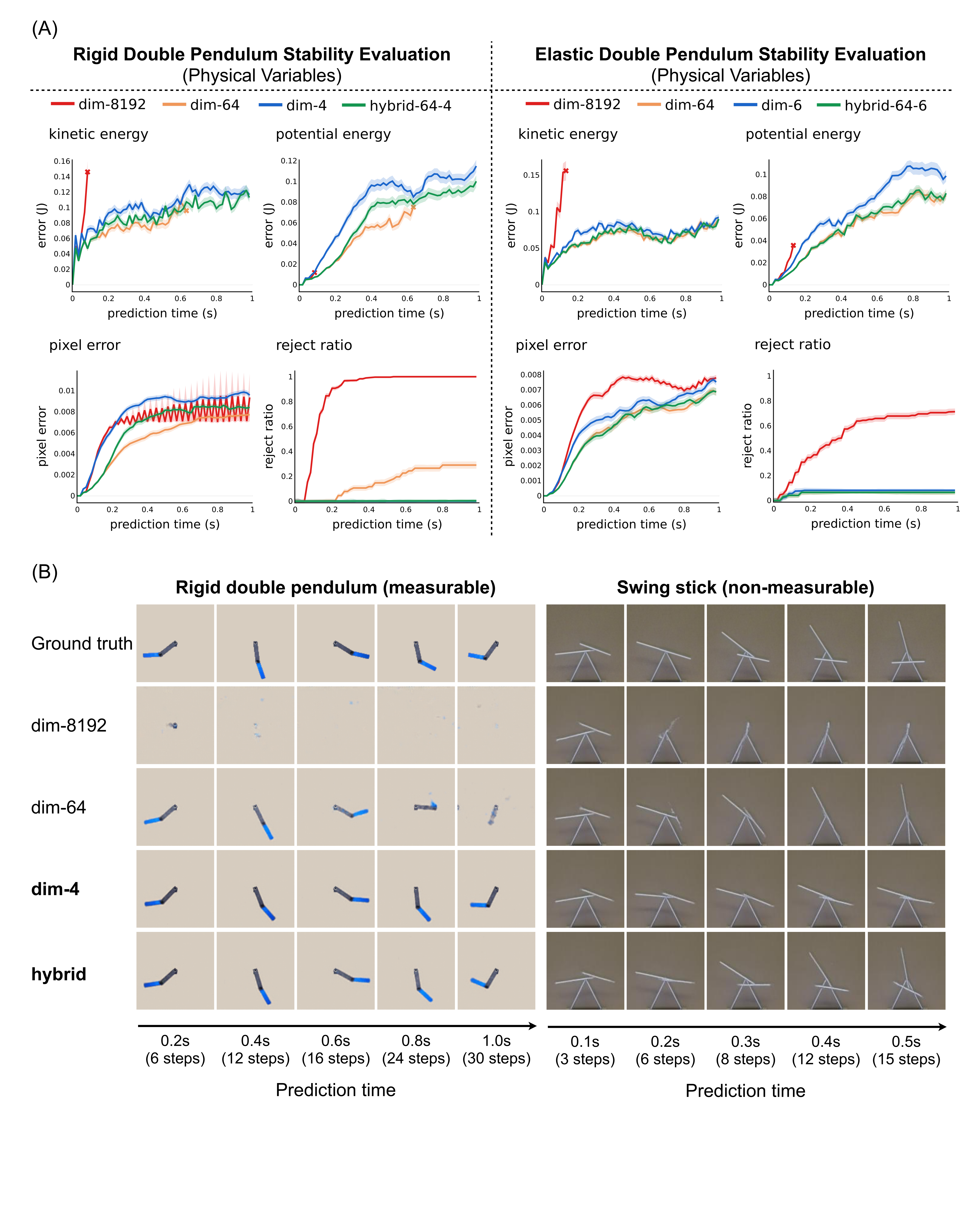}
    \caption{\textbf{Long-term Prediction Stability} (A) For systems where we could extract physical variables such as rigid and elastic double pendulum, using Neural State Variables as intermediate representations gives the most stable long-term future predictions. With our hybrid scheme, we can achieve 
    relatively stable and accurate predictions. Pixel error cannot replace physics-based evaluation for measuring the long-term prediction stability. (B) Predictions through Neural State Variables can maintain the ground truth state space for long-term predictions. The hybrid scheme offers 
    the most stable and accurate predictions.}
    \label{fig:reject-ratio}
    \vspace{-12pt}
\end{figure*}

Fig.~\ref{fig:reject-ratio}(A) shows the stability results on rigid double pendulum and elastic double pendulum where we can extract the physical variables from videos. The 8192-dim and 64-dim scheme cannot give stable long-term future predictions. In our experiments, we noticed that both schemes can provide stable predictions when the system intrinsic dimension is smaller or equal than 2. See Appendix for more results.

Inspired by lessons learned from computational physics, an effective fix to the unstable long-term prediction is to construct a prediction scheme where the predicted states will be projected into a small neighborhood of the state space. Here Neural State Variables serve as a strong candidate solution. This is because  Neural State Variables have the same dimension as the intrinsic dimension. This fact prevents predictions from falling off the system manifold into new dimensions. We provide more detailed theoretical analysis in the Appendix. 


The blue curves in Fig.~\ref{fig:reject-ratio}(A) illustrate the stability introduced by using Neural State Variables as intermediate representations for long-term predictions. Neural State Variables provide the most stable predictions across all the systems. However, since Neural State Variables were obtained by performing reconstruction on a relative high-dimensional latent embedding, they have an inferior performance on one-step prediction accuracy.

To combine the best of two worlds, we propose a hybrid scheme as our final solution: using Neural State Variables as stabilizers while performing long-term predictions with their corresponding high-dimensional latent embeddings. Formally, the hybrid scheme follows an $N + 1$ pattern where for every $N$ steps performed with the high-dimensional latent vectors, a one step prediction is followed with the Neural State Variables. As shown by the green curves in  Fig.~\ref{fig:reject-ratio}(A) and (B), our hybrid scheme offers 
stable and accurate long-term predictions. In Fig.~\ref{fig:reject-ratio}(A) and (B), the hybrid scheme was implemented with specific values of $N$ between 3 and 6. We also conducted experiments in the Appendix with different choices of $N$ and found that the outcomes were not sensitive to the particular values of $N$.

Another important note is that the use of pixel error as the evaluation metric, though easy to compute, can be misleading for the evaluation of long-term predictions when the predictions quickly become unstable. As quantitatively and qualitatively demonstrated in Fig.~\ref{fig:reject-ratio}, pixel errors remain roughly the same after the predicted images become plain backgrounds. These pixel errors are even smaller than the pixel errors computed from a slightly inaccurate but clear prediction. This observation further emphasizes the significance of designing an appropriate $M_{\mathcal{S}}(\cdot)$ metric.

\subsubsection*{Neural State Variables for Dynamics Stability Indicators}

So far, for evaluating long-term prediction stability, we have been assuming that we can extract the physical variables from the system states during evaluation. Yet, it is common that, in most of the video representations in our dataset, we know neither which variables to extract nor how to extract them directly from videos. As noted above, pixel errors are also not reliable. In this case, a very challenging but important problem is how we can evaluate the long-term prediction stability from videos.  Resolving this problem can potentially open up the door to quantitatively evaluate prediction stability of various schemes for many complex and unknown systems, all directly from videos.

Following the framework in the last section, the key is the design of the metric $M_{\mathcal{S}}(\cdot)$. Here we propose a solution based on Neural State Variables, namely $M_{\mathcal{S}}^{\textrm{neur}}$. Specifically, $M_{\mathcal{S}}^{\textrm{neur}}$ is a metric on a pair of states $(\hat{\vb*{X}}_t,\hat{\vb*{X}}_{t+dt})$.
\[ M_{\mathcal{S}}^{\textrm{neur}}(\hat{\vb*{X}}_t,\hat{\vb*{X}}_{t+dt})= \vecnorm{h_{\textrm{E}}\circ g_{\textrm{E}}(\hat{\vb*{X}}_{t+dt}) - \hat{F}_V(h_{\textrm{E}}\circ g_{\textrm{E}}(\hat{\vb*{X}}_t))},\]
where $\hat{F}_V$ is a neural network trained to approximate the latent dynamics on the space of Neural State Variables $\hat{\vb*{V}}_{t+dt\to t+2dt} \leftarrow \hat{F}_V (\vb*{V}_{t\to t+dt})$, $h_{\textrm{E}}\circ g_{\textrm{E}}(\hat{\vb*{X}}_t)=\hat{\vb*{V}}_{t\to t+dt}$ and $h_{\textrm{E}}\circ g_{\textrm{E}}(\hat{\vb*{X}}_{t+dt})=\hat{\vb*{V}}_{t+dt\to t+2dt}$ are Neural States in $\mathbb{R}^{\textrm{ID}}$, $\vecnorm{\cdot}$ is the Euclidean norm in $\mathbb{R}^{\textrm{ID}}$.

Intuitively, $M_{\mathcal{S}}^{\textrm{neur}}$ measures how far the predicted dynamics, reflected by the given predicted sequence, deviates from the reduced system dynamics projected onto the space of Neural State Variables. As shown in the previous sections, all latent embeddings with higher or equal to the intrinsic dimension may provide accurate short-term approximation of system dynamics. However, we chose the space of Neural State Variables as the reference because it has the same dimension as ID. First, as mentioned above, Neural State Variables project the predicted states in the small neighborhood of the ground truth states. Moreover, the Euclidean distance serve as a good metric to measure dynamics deviation in this case, while other higher dimensions may suffer the curse of dimensionality when designing the distance metric. Overall, $M_{\mathcal{S}}^{\textrm{neur}}$ is an ideal alternative candidate to $M_{\mathcal{S}}^{\textrm{phys}}$.

\begin{figure*}[t!]
    \centering
    \includegraphics[width=\linewidth]{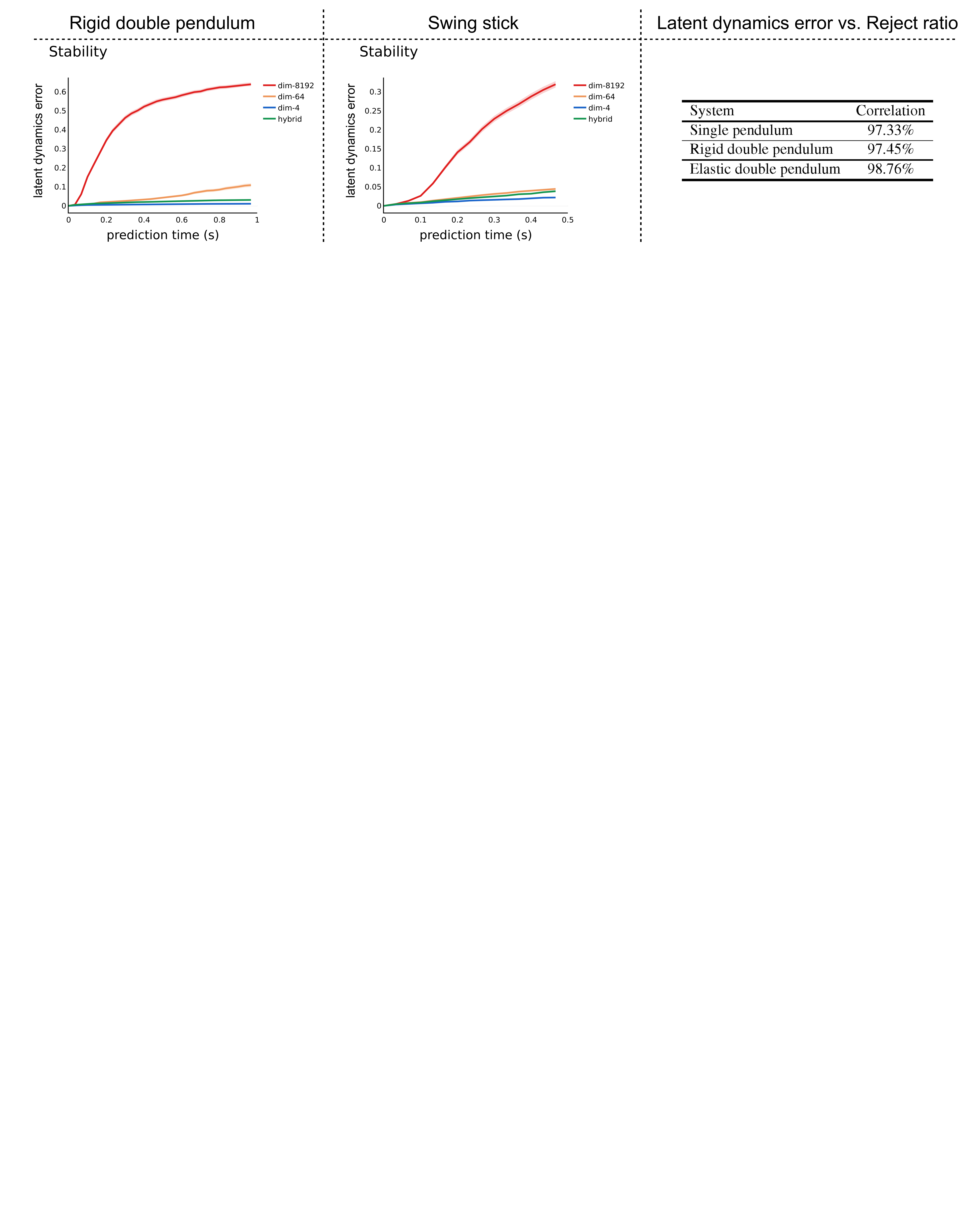}
    \caption{\textbf{Neural State Variables for Dynamics Stability Indicators} The latent dynamics on the space of Neural State Variables, namely neural latent dynamics, can be used as an effective indicator to quantify the stability of long-term future predictions. The strong correlation between the latent dynamics error and the physics reject ratio implies that our latent dynamics error can be used as an alternative metric to measure long-term prediction stability.}
    \label{fig:indicator}
\end{figure*}

Similar to $M_{\mathcal{S}}^{\textrm{phys}}$, the final $M_{\mathcal{S}}^{\textrm{neur}}$ is computed across multiple prediction sequences with various initial states. We show the evaluation results with our stability metric based on Neural State Variables in Fig.~\ref{fig:indicator}. $M_{\mathcal{S}}^{\textrm{neur}}$ produces patterns that highly match with $M_{\mathcal{S}}^{\textrm{phys}}$ for the systems where we know how to extract physical variables. This can also be seen in the correlation plot in Fig.~\ref{fig:indicator} where we computed the Pearson correlation coefficient between reject ratio of all models (dim-8192, dim-64, dim-ID, hybrid) at all prediction steps and the respective latent dynamics errors. For unknown systems, we observed the same trend where high-dimensional latent embeddings schemes are often not stable. In conclusion, our $M_{\mathcal{S}}^{\textrm{phys}}$ metric can help us measure the long-term prediction stability directly from videos without additional prior knowledge of the system.

\begin{figure*}[t!]
    \centering
    \includegraphics[width=0.97\linewidth]{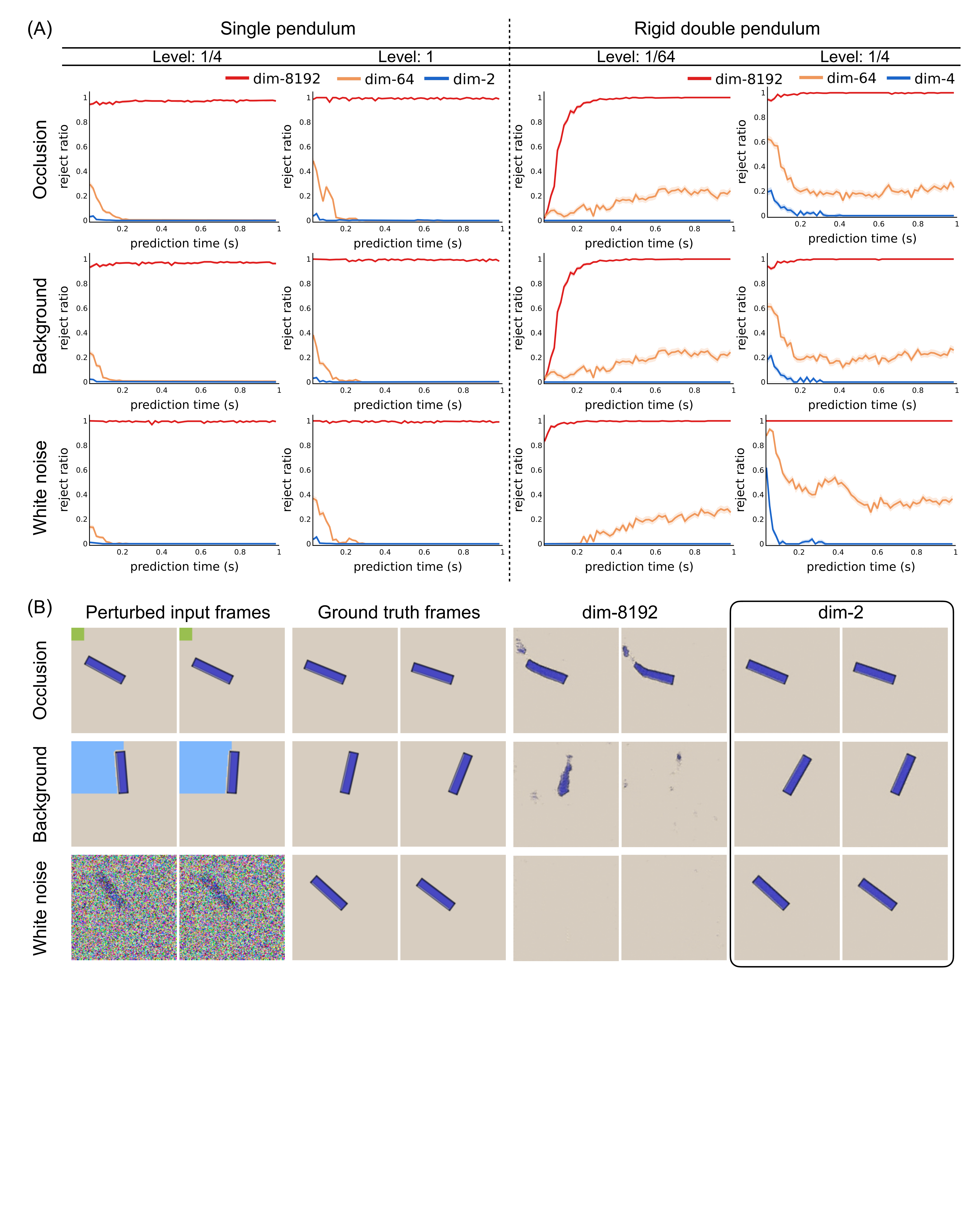}
    \caption{\textbf{Neural State Variables for robust long-term prediction} (A) The physical reject ratio can effectively measure long-term prediction robustness from perturbed initial data. For both systems and all combinations of perturbation types and perturbation levels, the model rollouts through Neural State Variables provide the most robust predictions. (B) From perturbed initial frames, the dynamics predictive model with high dimensional latent vectors produces blurred images or pure background with only one model rollout, while the model rollouts through Neural State Variables still give reasonable predictions.}
    \label{fig:robust}
\end{figure*}

\subsubsection*{Neural State Variables for Robust Long-Term Prediction}

Another critical factor when modeling system dynamics from videos is the robustness against visual perturbations. Therefore we applied several visual perturbations on the input video frames during test time and evaluated the performances of different models.

Specifically, we performed three types of perturbations. The first type is to simulate camera occlusions by covering certain portion of the input frames with a randomly generated color square. The area of the square indicates the level of the perturbation. For example, $\frac{1}{64}$ means the area of the square is $\frac{1}{64}$ times of the area of one input frame.

The second type is to simulate background color change by covering certain portion of the input frame background with a randomly generated color square. The main difference between this perturbation and the first one is that the color square will not cover the object. The level definition is the same with the first perturbation type.

Lastly, to simulate possible sensor noises, we added random Gaussian noise on the input frames. The Gaussian noise has a zero mean and different level of standard deviations. For example, $\frac{1}{64}$ means the Gaussian noise has a standard deviation of $\frac{1}{64} \times 255$ where $255$ is the highest pixel value in the input frames.

We show the test-time results using the physical reject ratio metric in Fig.~\ref{fig:robust}(A). The quantitative results clearly demonstrate the strong robustness of models on the Neural State Variables space across all level of perturbations. The models with very high dimensional latent space quickly produce unstable predictions.The models with a relative lower dimension but still higher than ID can sometimes give stable predictions again after several unstable rollouts.
However, even though the predictions can become stable again, it requires much more number of prediction steps. We also show qualitative visualizations in Fig.~\ref{fig:robust}(B).

\subsubsection*{Analysis}

We hypothesize that Neural State Variables contain rich physical meanings that align with the conventional definition of the physical State Variables. In this section, we verify this hypothesis through both quantitative regression experiments and qualitative visualizations.

\begin{figure*}[t!]
    \centering
    \includegraphics[width=\linewidth]{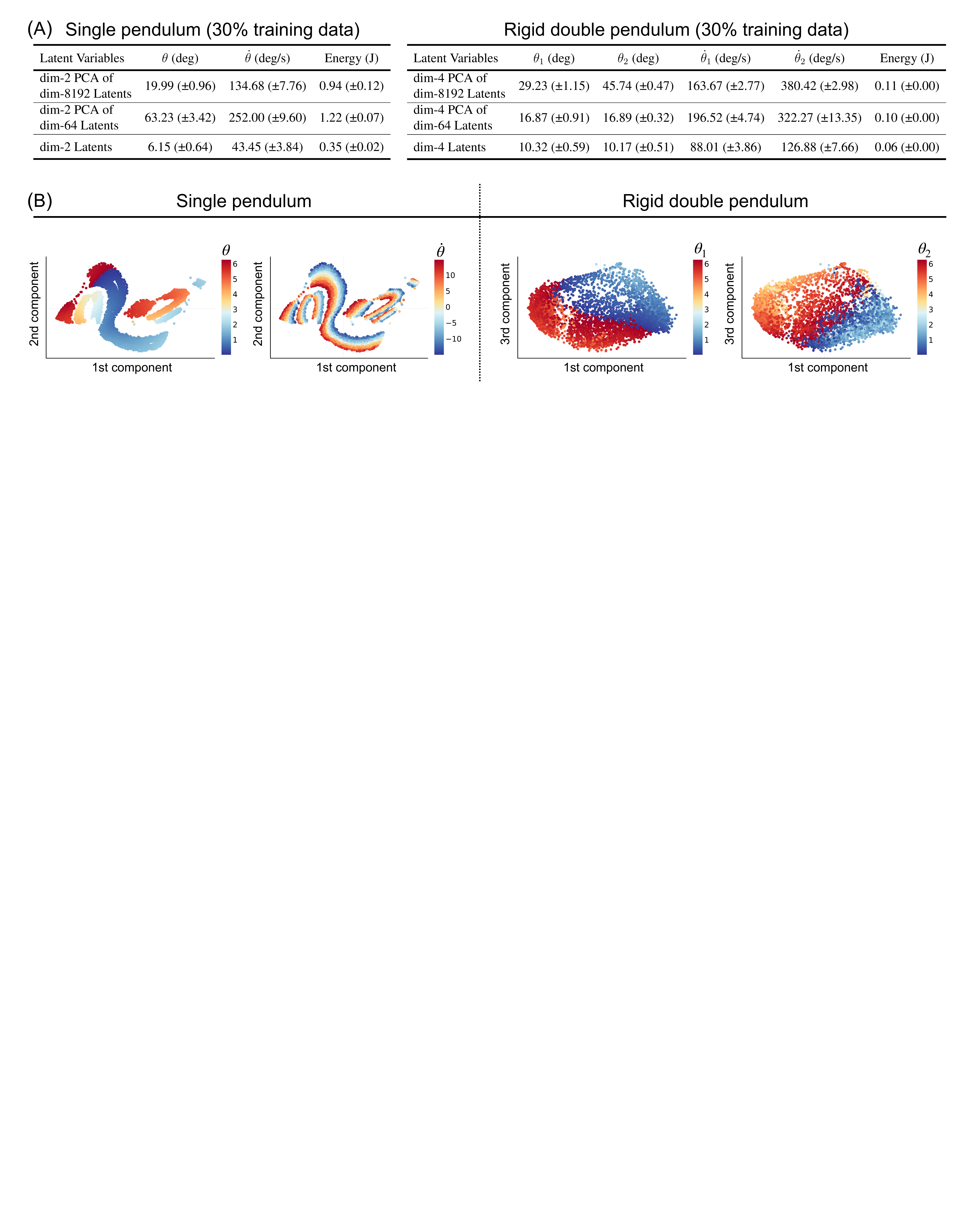}
    \caption{\textbf{PCA regression and Neural State Variables visualizations} (A) Neural State Variables capture much richer information about the system dynamics than state variables obtained through PCA from other high dimensional latent embedding vectors. (B) We also visualize the  Neural State Variables after applying PCA on them. The colors represent the value of different physical variables. Examples shown here suggest interesting symmetrical structures encoded in the Neural State Variables.}
    \label{fig:analysis}
\end{figure*}

We trained a small neural network with five layers of MLPs to regress conventional physical variables including positions, velocities, and energies from learned Neural State Variables. Our results are shown in Fig.~\ref{fig:analysis}(A). Using 30\% of labeled data, the learned Neural State Variables can be used to accurately regress the physical variables. We then compared the regression errors with those from the first few principal components of high dimensional latent vectors from our dynamics predictive model. Using the same number of state variables, which equals ID, and the same labeled data, the regression errors using principal components of high dimensional latent vectors are much larger than those using Neural State Variables, especially for velocity variables. Therefore, state variables obtained through PCA, or equivalently through linear neural networks, can hardly capture the dynamics of the system.

Visualizations colored by the value of physical variables in Fig.~\ref{fig:analysis}(B) can further demonstrate the physical meaning of the Neural State Variables. We observe that indeed the physical variables are captured in the set of Neural State Variables chosen by our modeling system. The charts also reveal the inherent symmetrical nature of these variables.

\subsubsection*{Outlook}

The burgeoning field of automated scientific discovery involves the development of machines that infer the underlying governing dynamics of various physical systems from observation data.
While the raw modeling capacity of machines clearly exceeds that of humans, the choice of what to model has remained human responsibility. We propose that the choice of state variables themselves may also be assisted or influenced by machines that can observe the system directly with data-rich sensors like cameras. Since machines can see the world in colors that we cannot see, hear the world in frequencies we cannot hear, and experience the world using senses we do not have, machines may help intuit new kinds of fundamental variables that we cannot imagine.

This work proposes several advancements that complement the existing works. First, it provides a way to distill a set of complete and non-redundant Neural State Variables from raw visual observations. Our work does not assume any prior knowledge about the physical system, in particular which variables to measure, or what mathematical primitives should be used for modeling. Second, this work shows the power of the Neural State Variables in producing stable long-term predictions against accumulation of model prediction errors as well as external visual perturbations.

There are many promising directions for the future research. First, the learned Neural State Variables can be further regularized with prior physical knowledge or existing physics-informed AI techniques, so that they could have better correspondence to traditional physical variables and satisfy physical constraints such as energy conservation \cite{pukrittayakamee2009simultaneous,wu2016physics,chmiela2017machine,schutt2017quantum,smith2017ani,bondesan2019learning,greydanus2019hamiltonian,lutter2019deep}.
Moreover, when frames are corrupted or contain incomplete information because of hidden factors, system uncertainty or inappropriate sampling frequency, the observation data may not fully capture the real physics. Regularization techniques \cite{lange2021fourier,mallen2021deep} can help us design algorithms to provide effective and automated remedies for such imperfect observation data. For complex systems with little prior physical knowledge, another interesting direction is to analyze the Neural State Variables and translate them into human interpretable physics. Finally, the proposed framework of distilling Neural State Variables and generating stable long term predictions can be used as a component of automated control systems.

\textbf{Funding:} This research was supported in part by NSF AI Institute for Dynamical Systems \#2112085, DARPA MTO Lifelong Learning Machines (L2M) Program HR0011-18-2-0020, NSF NRI \#1925157, NSF DMS \#1937254, NSF DMS \#2012562, NSF CCF \#1704833, and DE \#SC0022317. \textbf{Author contributions:} B.C. and H.L. proposed the research; B.C., K.H., H.L. and Q.D. performed experiments and numerical analysis, B.C. and K.H. designed the algorithms; B.C., K.H., I.C. and S.R. collected the dataset; B.C., K.H., H.L. and Q.D. wrote the paper; all authors provided feedback. \textbf{Data and materials availability:} We will open source all data and software codebase needed to evaluate the conclusion in the paper, the Supplementary Materials and our project website. Additional information can be addressed to B.C.

\bibliography{scibib}

\bibliographystyle{Science}

\clearpage

\end{document}